\documentclass[a4paper,12pt]{article}

\usepackage{amssymb}

\title{Relation of Origins of Primitive Chaos}

\author{Yoshihito Ogasawara\\
School of Fundamental Science and Engineering,\\
 Waseda University, Ohkubo, Shinjuku-ku, Tokyo 169-8555, Japan}

\begin{document}

\maketitle

{A new concept, primitive chaos, was proposed, as a concept closely related to the fundamental problems of sciences themselves such as determinism, causality, free will, predictability, and time asymmetry [{\em J. Phys. Soc. Jpn.} {\bf 2014}, {\em 83}, 1401]. This concept is literally a primitive chaos in such a sense that it leads to the characteristic properties of the conventional chaos under natural conditions. Then, two contrast concepts, nondegenerate Peano continuum and Cantor set, are known as the origins of the primitive chaos. In this study, the relation of these origins is investigated with the aid of a mathematical method, topology. Then, we can see the emergence of interesting concepts such as the relation of whole and part, and coarse graining, which imply the essence of our intrinsic recognition for phenomena.}

\section{Introduction}

The concept of primitive chaos was proposed, as a concept closely related to the fundamental problems of sciences themselves such as determinism, causality, free will, predictability, and time asymmetry \cite{O2010,O2014}.\\

{\noindent\bf Definition 1.}~~{\it If a set $X$, the family of subsets of $X$, $\{X_\lambda;~\emptyset\ne X_\lambda\subset X,~\lambda\in\Lambda\}$, and the family of maps, $\{f_{X_\lambda}:X_\lambda\to X,~\lambda\in\Lambda\}$, satisfy the following property (P), $(X,~\{X_\lambda,~\lambda\in\Lambda\},~\{f_{X_\lambda},~\lambda\in\Lambda\})$ is called primitive chaos.
\begin{itemize}
\item[(P)] For any infinite sequence $\omega_0,~\omega_1,~\omega_2,\ldots$, there exists an initial point $x_0\in\omega_0$ such that $f_{\omega_0}(x_0)\in\omega_1,~f_{\omega_1}(f_{\omega_0}(x_0))\in\omega_2,\ldots$, where $\omega_i\in\{X_\lambda,~\lambda\in\Lambda\}$ for each $i$.
\end{itemize}}~\\

In the primitive chaos, each set $X_\lambda$ represents an event or a selection, and each map $f_{X_\lambda}$ represents a law or causality. Then, under natural conditions, the primitive chaos leads to 
the characteristic properties of the conventional chaos, such as the existence of a nonperiodic orbit, the existence of the periodic point whose prime period is $n$ for any $n\in\mathbb{N}$, the existence of a dense orbit, the density of periodic points, sensitive dependence on initial conditions, and topological transitivity \cite{O2012,RLD}. In this sense, this primitive chaos is literally a primitive chaos.

Then, by exploring sufficient conditions for the guarantee of existence of the primitive chaos from a topological viewpoint, we saw the emergence of two contrast concepts, nondegenerate Peano continuum and Cantor set \cite{O2010,O2014}. In fact, if $X$ is a nondegenerate Peano continuum or a Cantor set, it guarantees not only the existence of the primitive chaos  $(X,~\{X_\lambda,~\lambda\in\Lambda\},~\{f_{X_\lambda},~\lambda\in\Lambda\})$, but also its infinite varieties. Namely, we can see such a situation that infinitely many events $X_\lambda,~\lambda\in \Lambda$ and laws $f_{X_\lambda},~\lambda\in \Lambda$ emerge from $X$, which satisfy the property (P). In this sense, we can state that the nondegenerate Peano continuum and the Cantor set are the origins of the primitive chaos.

Here, a nondegenerate Peano continuum is defined by a locally connected continuum more than one point, and a continuum is defined by a nonempty compact connected metric space \cite{nad}. Then, the nondegenerate Peano continuum is a universal concept, which has tremendous examples such as all arcs, all $n$-cells for all $n$, all $n$-spheres for all $n$, all toruses, all solid toruses, all trees, all graphs, all nondegenerate dendrites, all Hilbert cubes, and so on \cite{nad}. Therefore, taking it into consideration that the nondegenerate Peano continuum is the origin of the primitive chaos leading to the chaos, we can see an answer of a question ``Why are we surrounded by diverse chaotic behaviors?" \cite{O2010,O2014}

A Cantor set is defined by a topological space homeomorphic to the Cantor middle third set, and it is known that a topological space is a Cantor set if and only if it is a zero-dimensional perfect compact metrizable space \cite{nad}. Here, a topological space is perfect provided that it contains no isolated points, and a topological space is zero-dimensional provided that there is a base for its topology such that each element of the base is closed and open (clopen). Then, it is noted that the Cantor set is also a universal concept quite different from the special set, the Cantor middle-third set \cite{O2014}.

In addition, since the nondegenerate Peano continuum is characterized by its continuum and the Cantor set is characterized by its zero-dimensionality, the above results naturally remind us of two contrast aspects of matter from a macroscopic viewpoint and a microscopic viewpoint. More generally, the concepts of continuity (continuum) and discreteness (zero-dimensionality) seem to be our intrinsic concepts for the method of recognizing phenomena \cite{O2014}. 

In this study, we investigate the relation between these universal concepts, the nondegenerate Peano continuum and the Cantor set, which are the origins of the primitive chaos, by a mathematical method, topology. Here, it is noted that the topology is the rigorous method which can describe extremely abstract concepts, and it can be recognized as the generalization or abstraction of geometry. Therefore, we can interpret it as the discussion on concepts of morphology. However, in this study, we explore the possibility of interpreting it not only as such discussion, but also as the discussion on ``morphology of concepts" or as the discussion on morphology of our interior views \cite{CM,ZP,KH,KL,O2011,SE,PJ,TR}. Namely, presuming the symmetric relation between external and internal views, we are trying to overcome the dualistic relation of matter and mind.

\section{Relation of Origins of Primitive Chaos}

At first, let us recall the fact that if there exist finitely many maps more than two from the nondegenerate Peano continuum $X$ to itself such that the maps satisfy several conditions, there exists a Cantor set in $X$ \cite{ws}. Let us begin by generalizing this fact.\\

{\noindent\bf Proposition 2.}~~{\it Any nondegenerate Peano continuum $X$ has a Cantor set in it.}

{\noindent\bf Proof}~~By the definition, there exist two different points $x_1$ and $x_2$ in $X$. From Lemma 3 \cite{O2010,nad}, there exist two nondegenerate Peano subcontinua $X_1$ and $X_2$ of $X$ such that each $X_i$ contains $x_i$, and satisfies the relation
\begin{eqnarray}
dia\,X_i<\frac{d(x_1,x_2)}{3}\le\frac{dia\,X}{3}.
\end{eqnarray}
Let us represent
\begin{eqnarray}
X^1=X_1\cup X_2.
\end{eqnarray}
Then, for each $i$, there exist a point $x_{i2}$ in $X_i$, which is different from $x_i$, and let us represent
\begin{eqnarray}
x_{i1}=x_i.
\end{eqnarray}
From Lemma 3, there exist two nondegenerate Peano subcontinua $X_{i1}$ and $X_{i2}$ of $X_i$ such that each $X_{ij}$ contains $x_{ij}$, and satisfies the relation
\begin{eqnarray}
dia\,X_{ij}<\frac{d(x_{i1},x_{i2})}{3}\le\frac{dia\,X_i}{3}.
\end{eqnarray}
Let us represent
\begin{eqnarray}
X^2=X_{11}\cup X_{12}\cup X_{21}\cup X_{22}.
\end{eqnarray}
Then, this procedure can be repeated and a sequence $\{X^n\}$ is obtained. 

Let us verify that the set in $X$,
\begin{eqnarray}
S=\bigcap_nX^n,
\end{eqnarray}
is a Cantor set. First, since each $X^n$ is closed, $S$ is closed, and thus $S$ is compact. Then, for any point $x\in S$ and any $\varepsilon>0$, there exist a positive integer $n$ and $X_{i_1\cdots i_n}\subset X^n$ (each $i_j$ is 1 or 2) such that
\begin{eqnarray}
x\in X_{i_1\cdots i_n}\subset U(x,\varepsilon)=\{y\in X;~d(x,y)<\varepsilon\}.
\end{eqnarray}
Since there exist two different points $x_{i_1\cdots i_{n-1}1}$ and $x_{i_1\cdots i_{n-1}2}$ in $X_{i_1\cdots i_n}\cap S$, $S$ is perfect. In addition, $X_{i_1\cdots i_n}\cap S$ is an open subset of $S$ (i.e., it is a clopen subset of $S$), because 
\begin{eqnarray}
S-(X_{i_1\cdots i_n}\cap S)=(X^n-X_{i_1\cdots i_n})\cap S 
\end{eqnarray}
is a closed subset of $S$. Therefore, $S$ is zero-dimensional. $\Box$\\

{\noindent\bf Lemma 3.}~~{\it If $X$ is a nondegenerate Peano continuum, for any $\varepsilon>0$, there exist finitely many nondegenerate Peano subcontinua $X_1,\ldots,X_n$ of $X$ covering $X$ such that $dia\,X_i<\varepsilon$ for each $i$.}\\

Here, the nondegenerate Peano continuum has not only one Cantor set in it, but also infinitely many Cantor sets. In fact, the nondegenerate Peano subcontinuum $X_1$ of $X$ in the proof of Proposition 2 has a Cantor set $S_1$ in $X_1$, which is different from $S$. The nondegenerate Peano subcontinuum $X_{11}$ of $X$ also has a Cantor set $S_{11}$ in $X_{11}$, which is different from $S$ and $S_1$, and this discussion can be repeated. Accordingly, in addition to the fact that the nondegenerate Peano continuum has tremendous examples, each nondegenerate Peano continuum has infinitely many Cantor sets in it. Accordingly, we can see the relation of whole and part, as a relation between the Peano continuum and the Cantor set both of which are the origins of the primitive chaos. 

Furthermore, we can see such a situation that the whole is the coarse graining of the part. Let us first recall the concept of a decomposition space \cite{nad}.

{\noindent\bf Definition 4.}~~{\it For a topological space $(X,\tau)$ and a partition $\mathcal{D}$ of $X$, the family of subsets of $\mathcal{D}$,
\begin{eqnarray}
\tau(\mathcal{D})=\{\mathcal{U}\subset\mathcal{D};~\bigcup_{U\in\mathcal{U}}U\in\tau\},
\end{eqnarray}
is called a decomposition topology for $D$, and $(\mathcal{D},\tau(\mathcal{D}))$ is called a decomposition space of $X$.}\\

This concept is recognized as a coarse graining in such a sense that the decomposition topology is a natural topology for the partition from a topological viewpoint \cite{nad,O2011}. In addition, we can verify the following proposition.\\

{\noindent\bf Proposition 5.}~~{\it For a compact Hausdorff space $(X,\tau)$ and a partition  of $X$,
\begin{eqnarray}
\mathcal{D}=\{X_\lambda,~\lambda\in\Lambda\},
\end{eqnarray}
let a point $x_\lambda$ be the representative point of $X_\lambda$ for each $\lambda$. Then, the subspace of $X$ consisting of the representative points,
\begin{eqnarray}
Y=\{x_\lambda,~\lambda\in\Lambda\},
\end{eqnarray}
is homeomorphic to the decomposition space $(\mathcal{D},\tau(\mathcal{D}))$.}

{\noindent\bf Proof}~~The bijection
\begin{eqnarray}
h:Y\to\mathcal{D},\quad x_\lambda\to X_\lambda 
\end{eqnarray}
is continuous, because the relation
\begin{eqnarray}
h^{-1}(\mathcal{U})=\{x_\lambda;~X_\lambda\in \mathcal{U}\}=(\bigcup_{X_\lambda\in \mathcal{U}}X_\lambda)\cap Y
\end{eqnarray}
holds for any open subset $\mathcal{U}$ of $\mathcal{D}$. Since $X$ is compact and Hausdorff, $h$ is a homeomorphism. $\Box$\\

Here, we can recognize the subspace $Y$ of $X$ in Proposition 5 as a space coarsely grained by representing each set $X_\lambda$ as the point $x_\lambda$. By this proposition, we can see that the space $Y$ is the same as the decomposition space $\mathcal{D}$ from a topological viewpoint. In addition, it is noted that Proposition 2 holds regardless of the way of choice of each representative point $x_\lambda$ in $X_\lambda$.

Then, we can state that for any Cantor set $X$ and any nondegenerate Peano continuum $Y$, the decomposition space of $X$ is obtained, which is homeomorphic to $Y$. In fact, since there exists a continuous map $f$ from $X$ onto $Y$ from Lemma 6 \cite{O2014}, the decomposition space $\mathcal{D}_f$ is homeomorphic to $Y$ from Lemma 7 \cite{O2011}. Here, a topological space $X$ is a $T_1$-space provided that for each $x\in X$, the singleton $\{x\}$ is a closed subset of $X$, and thus any Hausdorff space is a $T_1$-space. \\

{\noindent\bf Lemma 6.}~~{\it Let $X$ be a zero-dimensional perfect 
compact $T_1$-space. For any compact metric space $Y$, there exists a continuous map from $X$ onto $Y$.}\\

{\noindent\bf Lemma 7.}~~{\it If $X$ is a compact space, $Y$ is a Hausdorff space, and $f$ is a continuous map $X$ onto $Y$, the decomposition space of $X$ for the partition $\{f^{-1}(y),~y\in Y\}~(=\mathcal{D}_f)$ is homeomorphic to $Y$.}\\

Accordingly, recalling that the relation of the nondegenerate Peano continuum and the Cantor set is the relation of whole and part, we can see such a picture that the coarsely grained part is identical to the whole from a topological viewpoint. While each nondegenerate Peano continuum has infinitely many Cantor sets in it, any Cantor set can be recognized as any nondegenerate Peano continuum by the concept of the decomposition space or the coarsely graining. 

Furthermore, let us see that there exist infinitely many continuous maps from any Cantor set onto any nondegenerate Peano continuum; that is, there exist infinitely many ways of such coarse graining.\\

{\noindent\bf Proposition 8.}~~{\it Let $X$ be a zero-dimensional perfect 
compact $T_1$-space and $Y$ be a compact metric space. For any positive integer $n$, any mutually disjoint nonempty clopen subsets $A_1,\ldots,A_n$ of $X$ covering $X$ and any nonempty closed subsets $B_1,\ldots,B_{n}$ of $Y$ covering $Y$, there exists a continuous map $f$ from $X$ onto $Y$ such that $f(A_i)=B_i$ for each $i$.}

{\noindent\bf Proof}~~Since the clopen subset of a zero-dimensional perfect 
compact space is also zero-dimensional, perfect, and 
compact, from Lemma 6, there exist continuous surjections
\begin{eqnarray}
g_i:A_i\to B_i,\quad i=1,\ldots,n.
\end{eqnarray}
Then, the surjection
\begin{eqnarray}
f:X\to Y,\quad x\mapsto
g_i(x),&x\in A_i
\end{eqnarray}
is continuous. In fact, since each $A_i$ is an open subset of $X$, for any open subset $U$ of $Y$,
\begin{eqnarray}
f^{-1}(U)=g_1^{-1}(U\cap B_1)\cup\cdots\cup g_n^{-1}(U\cap B_n)
\end{eqnarray}
is an open subset of $X$. $\Box$\\

Here, recalling the following lemma \cite{O2014} in addition to Lemma 3, we can see that there exist infinitely many continuous maps from the Cantor set onto the nondegenerate Peano continuum. Here, a $T_0$-space is a topological space $Y$ such that for any points $y_1$ and $y_2$ in $Y$ with $y_1\ne y_2$, there exists an open subset $U$ of $Y$ such that $y_i\in U$ and $y_j\notin U$ for some choice of $i$ and $j$, and thus any $T_1$-space is a $T_0$-space.\\

{\noindent\bf Lemma 9.}~~{\it Let $X$ be a zero-dimensional perfect $T_0$-space. For any positive integer $n$, there exist mutually disjoint nonempty clopen subsets $X_1,\ldots,X_n$ of $X$ covering $X$.}\\

In addition, Proposition 8 implies the following corollary because the complementary set of finite union of clopen sets is clopen.\\

{\noindent\bf Corollary 10.}~~{\it  If $X$ and $Y$ are Cantor sets, for any positive integer $n$, any $n$ mutually disjoint nonempty clopen subsets $A_1,\ldots,A_n$ of $X$, and any $n$ nonempty clopen  subsets $B_1,\ldots,B_n$ of $Y$, there exists a continuous map $f$ from $X$ onto $Y$ such that $f(A_i)=B_i$ for each $i$.}\\

This corollary guarantees the infinite varieties of the causality $\{f_{X_\lambda},~\lambda\in\Lambda\}$ in the primitive chaos $(X,~\{X_\lambda,~\lambda\in\Lambda\},~\{f_{X_\lambda},~\lambda\in\Lambda\})$ for the Cantor set $X$ \cite{O2014}. In fact, since a clopen subset of a Cantor set is again a Cantor set as its subspace, by Lemma 9, from one Cantor set $X$, the events or the selections $X_1,\ldots,X_n$ are generated with infinite varieties, which are Cantor sets as subspaces of $X$, and the causality $\{f_{X_1},\ldots,f_{X_n}\}$ are generated also with infinite varieties, by Corollary 10. 

In order to compare with the case of the nondegenerate Peano continuum, let us recall the following proposition \cite{O2010} which also guarantees the infinite varieties of the causality in the primitive chaos for the nondegenerate Peano continuum.\\

{\noindent\bf Proposition 11.}~~{\it If $X$ and $Y$ are nondegenerate Peano continua, for any positive integer $n$, any $n$ points $x_1,\ldots,x_n\in X$, and any $n$ points $y_1,\ldots,y_n\in Y$, there exists a continuous map $f$ from $X$ onto $Y$ such that $f(x_i)=y_i$ for each $i$.}\\

It is interesting that the law or the causality $f$ can be fitted to the ``points", i.e., $f(x_i)=y_i$, for the nondegenerate Peano continuum, while the law $f$ can only be fitted to the ``regions", i.e.,  $f(A_i)=B_i$, for the Cantor set. This relation again reminds us of two aspects of matter from a macroscopic viewpoint and a microscopic viewpoint, that is, the relation between the classical mechanics and the quantum mechanics.

\section{Conclusions}

The relation between the nondegenerate Peano continuum and the Cantor set is revealed, both of which guarantee the infinite varieties of the primitive chaos leading the chaos. Although they are contrast spaces (the one is characterized by continuum and the other is by zero-dimentionality), nevertheless they are closely related to each other. Any nondegenerate Peano continuum has infinitely many Cantor sets in it, and any Cantor set can be seen as any nondegenerate Peano continuum by the infinitely many ways of the coarsely graining. Furthermore, the result is obtained, which again reminds us of the aspects of matter from a macroscopic viewpoint and from a microscopic viewpoint. Recalling the fact that we are surrounded by diverse chaotic behaviors, these results not only indicate properties concerning the origin of the chaos, but also seem to imply the essence of our intrinsic recognition for phenomena.

\subsection*{Acknowledgments}
The author would like to acknowledge the support and useful comments of Professors Shin'ichi Oishi and Akira Koyama of Waseda University, and the helpful discussions of Professors Emeritus Yoshisuke Ueda of Kyoto University and Tsuneo Watanabe of Toho University. This study was supported by the Japan Science and Technology Agency.


\begin{thebibliography}{99}
\bibitem{O2010} Ogasawara, Y. Sufficient conditions for the existence
of a primitive chaotic behavior. {\em J. Phys. Soc. Jpn.} {\bf 2010}, {\em 79}, 15002.
\bibitem{O2014} Ogasawara, Y.; Oishi, S. Characteristic spaces emerging from primitive chaos. {\em J. Phys. Soc. Jpn.} {\bf 2014}, {\em 83}, 1401.
\bibitem{O2012} Ogasawara, Y.; Oishi, S. Consideration of a primitive chaos. {\em J. Phys. Soc. Jpn.} {\bf 2012}, {\em 81}, 103001.
\bibitem{RLD} Devaney, R.L. {\em An Introduction to Chaotic Dynamical Systems}; Westview Press: Colorado, 2003.

\bibitem{nad} Nadler, S.B.Jr. {\em Continuum Theory}; Marcel Dekker Inc: New York, 1992.

\bibitem{CM} Husserl, E.  {\em Cartesianische Meditationen}; Felix Meiner: Hamburg, 1992.
\bibitem{ZP} Husserl, E.  {\em Zur Phanomenologie des Inneren Zeitbewusstseins}; Martinus Nijhoff: Den Haag, 1966.
\bibitem{KH} Klaus, H.  {\em Lebendige Gegenwart}; Martinus Nijhoff: The Hague, 1966.
\bibitem{KL}Lewin, K.  {\em Principles of Topological Psychology}; McGraw-Hill: New York, 1936.
\bibitem{O2011} Ogasawara Y.; Oishi S. Addendum to eeSufficient conditions for the existence
of a primitive chaotic behaviorff. {\em J. Phys. Soc. Jpn.} {\bf 2011}, {\em 80}, 67002.
\bibitem{SE} Schrodinger, E. {\em Mind and Matter}; Cambridge University Press: Cambridge, U.K., 1958.
\bibitem{PJ} Piaget, J. {\em Le Structuralisme}; Presses Universitaires de France: Paris, 1968.
\bibitem{TR} Thom, R. {\em Stabilite Structurelle et Morphogenese}; InterEditions: Paris, 1977.


\bibitem{ws}  Kitada A.; Ogasawara Y. On a decomposition space of a weak self-similar set. {\em Chaos, Solitons and Fractals} {\bf 2005}, {\em 24}, 785-787; [Errata {\bf 2005}, {\em 25}, 1273].


\end{thebibliography}
\end{document}